\documentclass{amsart}
\usepackage{amssymb,graphpap}
\usepackage{graphics}
\usepackage{euler}
\usepackage{enumitem}
\newtheorem{theorem}{Theorem}[section]
\newtheorem{proposition}[theorem]{Proposition}

\theoremstyle{definition}

\theoremstyle{remark}

\newcommand{\U}{\mathbb{U}}
\newcommand{\C}{\mathbb{C}}

\newcommand{\Tr}{\operatorname{Tr}}
\newcommand{\Z}{\mathbb{Z}}

\newcommand{\Y}{\mathbb{Y}}

\newcommand{\CUE}{\operatorname{CUE}}
\newcommand{\M}{\mathbb{M}}
\newcommand{\B}{\mathbb{B}}
\newcommand{\W}{\mathcal{W}}
\newcommand{\R}{\mathbb{R}}
\newcommand{\E}{\mathbb{E}}
\newcommand{\SZ}{\operatorname{SZ}}

\begin{document}
\date{\today}
\sloppy

\title[Random Contractions]
{Vicious Walkers and Random Contraction Matrices}

\author[J. I. Novak]{Jonathan Novak}
\address{Queen's University, Department of Mathematics and Statistics,
Jeffery Hall, Kingston, ON K7L 3N6, Canada}
\email{jnovak@mast.queensu.ca}
\urladdr{www.mast.queensu.ca/$\sim$jnovak}
%\subjclass[2000]{46L54 (46L65, 46L53, 60G09)}
%\keywords{Free probability, quantum exchangea
	
	\begin{abstract}
		The ensemble $\CUE^{(q)}$ of truncated random unitary matrices is a deformation
		of the usual Circular Unitary Ensemble depending on a discrete non-negative parameter
		$q.$  $\CUE^{(q)}$ is an exactly solved model of random contraction matrices originally
		introduced in the context of scattering theory.
		In this article, we exhibit a connection between $\CUE^{(q)}$ and Fisher's random-turns
		vicious walker model from statistical mechanics.  In particular, we show that 
		the moment generating function of the trace of a random matrix from $\CUE^{(q)}$ 
		is a generating series for the partition function of Fisher's model, when the 
		walkers are assumed to represent mutually attracting particles.
	\end{abstract}

	\maketitle
	
	\section{Introduction}

\subsection{Truncated random unitary matrices}
Fix an integer $d \geq 1,$ and let $\M_d=\M_d(\C)$ be the space of $d \times d$ complex matrices.
Consider the linear map
\begin{equation}
	T:\M_{d+1} \rightarrow \M_d
\end{equation}
which acts by removing the last row and column of a matrix.
For example, $T:\M_3 \rightarrow \M_2$ is defined by
\begin{equation}
	T\bigg{(} \begin{bmatrix} z_{11} & z_{12} & z_{13} \\
			z_{21} & z_{22} & z_{23} \\
			z_{31} & z_{32} & z_{33} \end{bmatrix} \bigg{)} =
			\begin{bmatrix} z_{11} & z_{12} \\
				z_{21} & z_{22}
			\end{bmatrix}
\end{equation}
$T$ is called the {\it truncation} map.  Since $\|T(M)\| \leq \|M\|$ in operator norm
for any $M \in \M_{d+1},$ $T$ maps the unitary group
\begin{equation}
	\U_{d+1}=\{U \in \M_{d+1} : U^*=U^{-1} \}
\end{equation}
into  
\begin{equation}
	\B_d=\{P \in \M_d : \|P\| \leq 1\},
\end{equation}
the semigroup of linear contractions of Euclidean space $\C^d.$
More generally, for any integer $q \geq 0,$ the $q$-fold composition $T^{(q)}$
maps each matrix in $\M_{d+q}$ to its $d \times d$ principal submatrix and thus sends
$\U_{d+q}$ into $\B_d.$

The truncation operator induces a very natural deformation of the Circular Unitary 
Ensemble ($\CUE$) from random matrix theory.  Consider the unitary group
$\U_{d+q}$ as a probability space, with the Borel $\sigma$-algebra and 
Haar probability measure.  The pushforward
$\gamma_d^{(q)}$ of Haar measure on $\U_{d+q}$ under $T^{(q)}$
is a Borel probability measure on $\B_d,$ and thus one obtains a random matrix 
ensemble
\begin{equation}
	\{(\B_d, \gamma_d^{(q)}) : d \geq 1\}
\end{equation}
which will be denoted $\CUE^{(q)}$ and called an ensemble of {\it truncated random
unitary matrices}.  Since $T^{(0)}$ is the identity operator on $\M_d,$ $\CUE^{(0)}$
reduces to the usual Circular Unitary Ensemble of unitary matrices under Haar measure
when $q=0.$

Ensembles of truncated random unitary matrices were first studied by Sommers and 
Zyczkowski \cite{s-z} in the context of quantum chaotic scattering, and
have been further investigated
by Petz and R\'effy \cite{p-r1}, \cite{p-r2}, Fyodorov and Khoruzhenko \cite{f-k}, and 
Neretin \cite{neretin}.  They are an important technical ingredient in Krishnapur's recent
study of random matrix-valued analytic functions \cite{krishnapur}.

Certain averages over $\CUE$ are well-known to be related to the combinatorics of
increasing subsequences in permutations \cite{rains}, \cite{b-r} and vicious walkers on $\Z$
\cite{a-vm}, \cite{f-g}.  A natural question is whether (and in what form) this relationship
with combinatorics extends to the deformed ensemble $\CUE^{(q)}.$

\subsection{Fisher's random-turns model}
\label{RTmodel}
Vicious walker models were introduced in statistical mechanics by Fisher in order to
model wetting and melting \cite{fisher}.
Fisher's {\it random-turns} vicious walker model consists of a system of $d$ particles
(``walkers'') initially occupying sites
\begin{equation}
	\mu=\mu_1>\mu_2>\dots>\mu_d
\end{equation}
on the integer lattice $\Z.$  Note that the walkers are labelled from right to left, so that
we list their positions in decreasing order.  The system evolves in discrete time according to the following
rule: at each instant, a single random particle makes
a random unit jump left or right (a ``random turn''), subject only to the constraint that no two particles can 
occupy the same lattice site simultaneously.  
The function $Z_d(N;\mu,\lambda)$ which counts the number of ways in which $d$ random-turns particles 
can depart initial sites $\mu$ and arrive at new sites
\begin{equation}
	\lambda=\lambda_1>\lambda_2>\dots>\lambda_d
\end{equation}
at time $N$ is the {\it partition function} of the model.  

If one assumes the existence of an 
attractive force between the particles, then the system can be at equilibrium only when 
the particles are on adjacent sites.  In this case we use the simplified notation 
$Z_d(N;q)$ to denote the number of ways in which $d$ random-turns particles can move
between ground states $q$ sites apart in $N$ instants.  Figure 1 gives an
example of a sequence of configurations of mutually attracting random-turns particles
counted by $Z_3(10;2).$

\begin{figure}
	\begin{center}
	\scalebox{0.5}[0.5]{
	\includegraphics{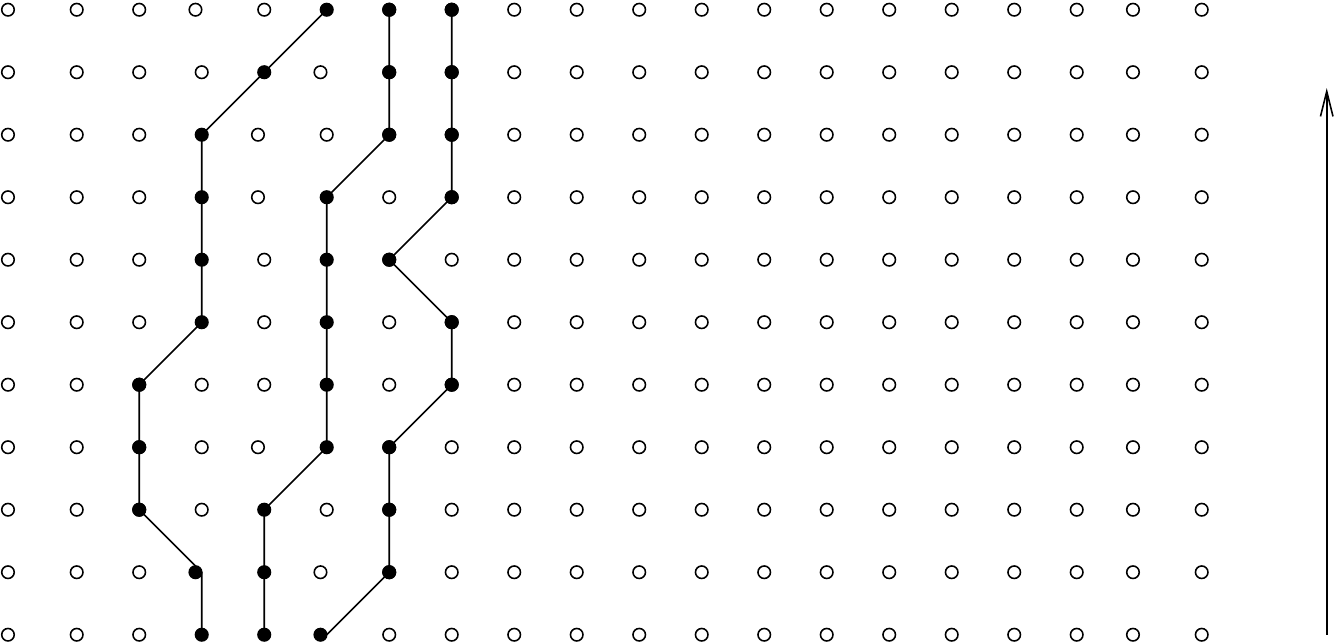}
	}
	\caption{\it Three random-turns particles moving between ground states.}
	\end{center}
\end{figure}

\subsection{Main result}
Our main result is the following.  Consider the matrix integral
\begin{equation}
	\label{integral}
	G_d(x;q):=\frac{x^{dq}}{H_{d \times q}} \int_{\B_d} e^{x\Tr(P+P^*)} \gamma_d^{(q)}(dP),
\end{equation}
where
\begin{equation}
	H_{d \times q}=\prod_{i=0}^{d-1} \frac{(q+i)!}{i!}
\end{equation}
is the hook-product of the $d \times q$ rectangular Young diagram.  $G_d(x;q)$ is 
a scaled and shifted version of the moment generating function of $\Tr(P_d^{(q)}+P_d^{(q)*}),$ 
where $P_d^{(q)}$ is a random matrix from $\CUE^{(q)}.$  

\begin{theorem}
	\label{main}
	$G_d(x;q)$ is the exponential generating series for the partition function of a system of $d$
	mutually attracting random-turns particles:
	\begin{equation}
		G_d(x;q)=\sum_{N \geq 0} Z_d(N;q) \frac{x^N}{N!}.
	\end{equation}
\end{theorem}

Below, we use this interpretation of $G_d(n;q)$ as a generating series to give a completely combinatorial proof of the following
identity, which expresses the matrix integral $G_d(x;q)$ as a Toeplitz determinant of 
Bessel functions.

\begin{theorem}
	\label{secondary}
	We have
	\begin{equation}
		G_d(x;q)=\det(I_{q+j-i}(2x))_{1 \leq i,j \leq d},
	\end{equation}
	where the entries in the determinant are modified Bessel functions.
\end{theorem}

{\it Remark:} Neretin \cite{neretin} has shown that, when $q \geq d, $ the measure 
$\gamma_d^{(q)}$ becomes absolutely continuous with respect to Lebesgue measure
on $\B_d$ and has density
\begin{equation}
	\frac{H_{d\times q}}{\pi^{d^2}H_{d \times (q-d)}} \det(I-P^*P)^{q-d} dP.
\end{equation}
Therefore when $q \geq d,$ $G_d(x;q)$ can be written more concretely as
\begin{equation}
	G_d(x;q)=\frac{x^{dq}}{\pi^{d^2}H_{d \times (q-d)}} \int_{\B_d} 
	e^{x\Tr(P+P^*)} \det(I-P^*P)^{q-d}dP.
\end{equation}
		
\subsection{Connection with the increasing subsequence problem}
Recall that a permutation $\sigma$ from the symmetric group $S(n)$ is said to 
have an {\it increasing subsequence} of length $k$ if there exist indices
\begin{equation}
	1 \leq i_1<i_2< \dots < i_k \leq n
\end{equation}
such that
\begin{equation}
	\sigma(i_1) < \sigma(i_2) < \dots < \sigma(i_k).
\end{equation}
Increasing subsequences in permutations were first studied by Erd\"os and Szekeres
\cite{e-s} in the $1930$s in connection with a Ramsey-type problem for the 
permutation group.  In the $1960$s,
Ulam \cite{ulam} raised the problem of determining the number $u_d(n)$ of permutations
in $S(n)$ with increasing subsequence length bounded by $d.$  This came to be known
as the {\it increasing subsequence problem}.  Stanley's ICM contribution \cite{stanley:icm} 
gives a comprehensive survey, from a combinatorial perspective, of the vast literature on increasing subsequences in permutations and related topics.  

Forrester \cite{forrester} has observed the following connection between the increasing 
subsequence problem and Fisher's random-turns model:
\begin{equation}
	\label{Forrester}
	Z_d(N;0)=\begin{cases} 
		{2n \choose n}u_d(n), \text{ if $N=2n$ for some $n \geq 0$}\\
		0, \text{ otherwise}
	\end{cases}.
\end{equation}

We will present a new proof of (a more general version of) Forrester's result below.  
Our approach 
is based on the fact that the configuration space of Fisher's random-turns model naturally
carries the structure of a graded graph with commutative raising and lowering operators.
This approach was inspired by Stanley's notion of a ``differential poset,'' defined as a 
graded graph which arises as the Hasse graph of a poset and whose raising and lowering operators satisfy the Heisenberg commutation
relation \cite{stanley:dp}.
In particular, this observation implies that specializing $q=0$ in Theorems \ref{main}
and \ref{secondary}, we obtain
\begin{equation}
	G_d(x;0)=\sum_{n \geq 0} u_d(n) \frac{x^{2n}}{n!n!}=det(I_{j-i}(2x))_{1 \leq i,j \leq d}.
\end{equation}
The identity between the matrix integral and the generating series of $u_d(n)$ is due to 
Rains \cite{rains}, while the identity between the series and the determinant 
is due to Gessel \cite{gessel}.  Both were discovered in the context of the enumeration 
of permutations with bounded increasing subsequence length.

\section{Configuration Space as a Graded Graph}
\label{config}

\subsection{Configuration space of random-turns particles}
The configuration space of a physical system is the set of all possible positions of its
constituents.  In the case of the random-turns model, this is the discrete space
\begin{equation}
	\W_d=\{(\lambda_1,\lambda_2,\dots,\lambda_d) \in \Z^d : \lambda_1>\lambda_2>
	\dots >\lambda_d\}.
\end{equation}
$\W_d$ is familiar as a type $A$ {\it Weyl lattice}, i.e. it is the intersection of 
$\Z^d$ with an open Weyl chamber for the type $A$ root system in $\R^d$ 
(see e.g. \cite{bump} about root systems and their Weyl chambers).
$\W_d$ becomes a simple, connected, locally finite graph when we declare
vertices $\mu,\lambda \in \W_d$ adjacent if and only if they are unit Euclidean distance
apart.  The partition function $Z_d(N;\mu,\lambda)$ of the random-turns model
counts the number of walks of length $N$ from $\mu$ to $\lambda$
on the graph $\W_d.$

$\W_d$ is also an example of a {\it graded graph}.
For our purposes, it suffices to define a graded graph to be a pair $(G,r)$ consisting
of a simple, connected, locally finite graph $G$ together with a function 
$r:G \rightarrow \Z,$ called the {\it rank function}, which associates an integer to each 
vertex of $G.$  The incidence relation and the rank function on $G$ are required to be 
compatible in the sense that if $u,v \in G$ are adjacent vertices, then either 
$r(v)=r(u)+1$ (denoted $u \nearrow v$) or $r(v)=r(u)-1$ (denoted $u \searrow v$).
Graded graphs occur frequently in combinatorics as the Hasse graphs of posets 
(see \cite{stanley:dp} and \cite{fomin}) and in representation theory where they are
known as {\it branching graphs} or {\it Bratelli diagrams} (see e.g. \cite{k-o-o}).

If we define a rank function $r$ on $\W_d$ by
\begin{equation}
	r(\lambda)=r(\lambda_1,\dots,\lambda_d):=\sum_{i=1}^d \lambda_i - {d+1 \choose 2},
\end{equation}
then $(\W_d,r)$ is a graded graph.  Note that we have normalized the rank function in 
such a way that the rank of the canonical ground state $\rho=(d,d-1,\dots,1)$ 
in Fisher's model is $0.$

Let $(G,r)$ be a graded graph.  By the {\it unrefined} partition function of $G$
we mean the number $Z_G(N;u,v)$ of walks on $G$ from $u$ to $v.$  We will refer to the 
number
$Z_G(L^{b_k}R^{a_k} \dots L^{b_1}R^{a_1};u,v)$ of
walks on $G$ from $u$ to $v$ of the form
\begin{equation}
	u \underbrace{\nearrow \nearrow \dots \nearrow}_{a_1}
	\underbrace{\searrow \searrow \dots \searrow}_{b_1}
	\dots
	\underbrace{\nearrow \nearrow \dots \nearrow}_{a_k}
	\underbrace{\searrow \searrow \dots \searrow}_{b_k} v
\end{equation}
as the {\it refined} partition function of $G.$
In the special case $G=\W_d$ we will continue to use the notation
$Z_{d}(N;\mu,\lambda)=Z_{\W_d}(N;\mu,\lambda)$ from the Introduction, as well
as $Z_d(L^{b_k}R^{a_k} \dots L^{b_1}R^{a_1};\mu,\lambda):=
Z_{\W_d}(L^{b_k}R^{a_k} \dots L^{b_1}R^{a_1};\mu,\lambda).$

$L$ and $R$ as defined above are simply inert symbols recording instances of 
positive and negative steps on the configuration graph (equivalently left 
and right jumps of the particles).  We may represent $L$ and $R$ as linear operators
acting on $\C[G],$ the free $\C$-vector space spanned by the vertices of $G,$ 
by defining
\begin{equation}
	R(u) := \sum_{u \nearrow v} v \text{ and } L(u):= \sum_{u \searrow w}w
\end{equation}
and extending linearly over $\C[G].$
$R$ and $L$ are then called the {\it raising} and {\it lowering} operators on $G.$
The monoid $\{L,R\}^*$ generated by the raising and lowering operators determines
the refined partition function, since
\begin{equation}
	Z_G(L^{b_k}R^{a_k} \dots L^{b_1}R^{a_1};u,v)=
	[v]L^{b_k}R^{a_k} \dots L^{b_1}R^{a_1}(u)
\end{equation}
where $[v]$ is the ``coefficient of $v$'' functional on $\C[G].$

\begin{theorem}
	\label{commuting}
	The raising and lowering operators on $(\W_d,r)$ commute.
\end{theorem}

\begin{proof}
	Let $\overline{\W}_d$ be the graded graph with vertex set
	\begin{equation}
		\overline{\W}_d=\{(\lambda_1,\lambda_2,\dots,\lambda_d) \in \Z^d :
		\lambda_1 \geq \lambda_2 \geq \dots \geq \lambda_d\},
	\end{equation}
	incidence relation as with $\W_d,$ and rank function defined by
	\begin{equation}
		\overline{r}(\lambda)=\overline{r}(\lambda_1,\lambda_2,\dots,\lambda_d):=
		\sum_{i=1}^r \lambda_i.
	\end{equation}
	Since the translation
		\begin{equation}
			\lambda \mapsto \overline{\lambda}:=\lambda-\rho
		\end{equation}
	is a rank-preserving isomorphism of $\W_d$ with $\overline{\W}_d,$ it suffices 
	to prove that the raising and lowering operators associated to 
	$(\overline{W}_d,\overline{r})$ commute.
	
	Let $\mu,\lambda \in \overline{\W}_d$ be arbitrary vertices.  Consider separately the 
	cases $\mu \neq \lambda$ and $\mu = \lambda.$
	
	\begin{description}
		
		\item[\bf Case $\mu \neq \lambda$]
		$[\lambda]LR(\mu)$ counts the number of walks on $\W_d$ from $\mu$ to 
		$\lambda$ of the form
		\begin{equation}
			\mu \nearrow \mu+e_i \searrow \mu + e_i -e_j = \lambda,
		\end{equation}		
		where $e_i$ and $e_j$ are standard basis vectors of $\R^n.$  If such
		an $i$ and $j$ exist, then $i \neq j$ by the assumption $\mu \neq \lambda.$  Thus
		\begin{equation}
			\mu \searrow \mu-e_j \nearrow \mu-e_j+e_i
		\end{equation}
		is a valid walk on $\W_d$ from $\mu$ to $\lambda,$ and we have a bijection
		\begin{equation}
			\mu \nearrow \nu \searrow \lambda \iff 
			\mu \searrow \nu' \nearrow \lambda
		\end{equation}
		between up-down walks on $\W_d$ from $\mu$ to $\lambda$ and 
		down-up walks on $\W_d$ from $\mu$ to $\lambda.$  Hence
		$[\lambda]LR(\mu)=[\lambda]RL(\mu).$
		
		\item[\bf Case $\mu = \lambda$]  In this case, let $r$ be the number of 
		distinct entries of $\mu=\lambda,$ i.e. $i_1<i_2<\dots <i_r$ satisfy
		$\mu_{i_1}>\mu_{i_2}> \dots >\mu_{i_r}.$  Then it is clear that
		$[\mu]RL(\mu)=[\mu]LR(\mu)=r.$
	\end{description}
\end{proof}

\subsection{Determinant identities}
Using Theorem \ref{commuting}, we can easily prove determinantal identities for the 
refined partition function of random-turns particles.
Let $I_k$ be the modified Bessel function of order $k,$
\begin{equation}
	\label{Bessel}
	I_k(2x) := \sum_{n \geq 0} \frac{x^n}{\Gamma(n+1)} \frac{x^{n+k}}{\Gamma(n+k+1)}.
\end{equation}
It is known that $I_{-k}(2x)=I_k(2x)$ for $k \in \Z$ (see e.g. \cite{b-m-m-s} about Bessel 
functions and their properties). 

\begin{proposition}
	\label{determinantal}
Let $\mu,\lambda \in \W_d$ be configurations with $r(\mu) \leq r(\lambda),$ and
let $W_0,W_1,\dots,W_n,\dots \in \{L,R\}^*$ be any sequence of words in the raising 
and lowering operators on $\W_d$ verifying 
\begin{equation}
	\deg_L W_n = n \text{ and } \deg_R W_n = n+r(\lambda)-r(\mu).
\end{equation}
Then
\begin{equation}
	\sum_{n \geq 0} Z_d(W_n;\mu,\lambda) \frac{x^{2n+r(\lambda)-r(\mu)}}
	{n!(n+r(\lambda)-r(\mu))!} = \det (I_{\lambda_i-\mu_j}(2x))_{1 \leq i,j \leq d}.
\end{equation}
\end{proposition}

\begin{proof}
	Since $\W_d$ is the intersection of $\Z^d$ 
	with an open type $A$ Weyl chamber in $\R^d,$ it follows immediately from the 
	Andr\'e-Gessel-Zeilberger reflection principle \cite{g-z}, \cite{g-m}, \cite{g-w-w} 
	that the generating series for the unrefined partition function is
		\begin{equation}
			\sum_{N \geq 0} Z_d(N;\mu,\lambda) \frac{x^N}{N!}=
			\det(I_{\lambda_i-\mu_j}(2x))_{1 \leq i,j \leq d}.
		\end{equation}
	Now, an $N$-step walk from $\mu$ to $\lambda$ on $\W_d$ exists if and only if 
	$N=2n+r(\lambda)-r(\mu)$ for some $n \geq 0$ (this $n$ being the number of negative
	steps).  By Theorem \ref{commuting}, the number of such walks is
		\begin{equation}
			{2n + r(\lambda)-r(\mu) \choose n}Z_d(W_n;\mu,\lambda)
		\end{equation}
	for any $W_n \in \{L,R\}^*$ with $\deg_L W_n=n$ and $\deg_R W_n=n+r(\lambda)-r(\mu).$
	Thus
	\begin{align}
		\sum_{N \geq 0} Z_d(N;\mu,\lambda) \frac{x^N}{N!}&=
		\sum_{n \geq 0}{2n + r(\lambda)-r(\mu) \choose n}Z_d(W_n;\mu,\lambda)
		\frac{x^{2n+r(\lambda)-r(\mu)}}{(2n+r(\lambda)-r(\mu))!}\\
		&=\sum_{n \geq 0} Z_d(W_n;\mu,\lambda)\frac{x^{2n+r(\lambda)-r(\mu)}}{n!
		(n+r(\lambda)-r(\mu))!}
	\end{align}
	and the result follows.
\end{proof}

Corollary \ref{determinantal} generalizes recent results of Xin (\cite{xin}, Theorem 11 and 
Proposition 12), who proves the case $W_n=L^nR^n$ using the ``Stanton-Stembridge trick.''

\subsection{Young tableaux and increasing subsequences}
Let $\Y$ be the {\it Young graph}, i.e. the Hasse graph of the lattice of Young diagrams partially
ordered by inclusion of diagrams (see \cite{stanley:ec2}).  $\Y$ is a graded graph, where the
rank $|\lambda|$ of a Young diagram $\lambda \in \Y$ is the number of cells in $\lambda.$
$\Y$ has a unique vertex $\emptyset$ of rank $0$ (the ``empty diagram''), and all other vertices
have positive rank.  Walks on the Young graph are known as {\it oscillating Young tableaux}.

Let $\Y_d$ be the induced subgraph of $\Y$ whose vertices are the Young diagrams with at 
most $d$ rows.  It is a well-known consequence of the RSK correspondence
(see \cite{stanley:ec2} as well as the original article of Schensted \cite{schensted}) that
\begin{equation}
	\label{RSK}
	Z_{\Y_d}(L^nR^n;\emptyset,\emptyset)=u_d(n),
\end{equation}
where $u_d(n)$ is the number of permutations in $S(n)$ with no increasing subsequence
of length greater than $d.$
		
Observe that there is a canonical embedding $\imath: \Y_d \hookrightarrow \overline{W}_d$
obtained by mapping each Young diagram $\lambda$ in $\Y_d$ onto the vector of 
its row lengths.  Composing $\imath$ with the translation
\begin{equation}
	\lambda \mapsto \lambda^{\circ} := \lambda+\rho
\end{equation}
gives an embedding of $\Y_d$ into $\W_d.$  Thus it is an immediate consequence of 
Theorem \ref{commuting} that
\begin{equation}
	Z_d(N;0)=\begin{cases}
		{2n \choose n}u_d(n), \text{ if $N=2n$}\\
		0, \text{ otherwise}
	\end{cases}
\end{equation}
which is precisely Forrester's result (\ref{Forrester}).  We therefore have the following Proposition
as a direct consequence of Proposition \ref{determinantal}.

\begin{proposition}
	\label{Gessel}
	For any $d \geq 1,$ 
	\begin{equation}
		\sum_{n \geq 0} u_d(n) \frac{x^{2n}}{n!n!} = \det(I_{i-j}(2x))_{1 \leq i,j \leq d}.
	\end{equation}
\end{proposition}

\begin{proof}
	Choose $\mu=\lambda=\rho$ in Corollary \ref{determinantal}.
\end{proof}

This is precisely Gessel's identity, the original proof of which appears in \cite{gessel}.
Alternative proof of this result were later given by
Gessel, Weinstein, and Wilf \cite{w-w}, Tracy and Widom \cite{t-w},
and Xin \cite{xin}.  Gessel's identity was the starting point of Baik, Deift, and Johansson \cite{b-d-j} in their groundbreaking work on the limiting distribution of the length of the longest increasing subsequence
in a large random permutation.

\section{Proof of the main theorem}
\label{proofMain}

We now give the proof of Theorem \ref{main}, which links truncated random unitary 
matrices with Fisher's random-turns model.  The proof is based on the integral identity
\begin{equation}
	\label{identity}
	G_d(x;q)=\int_{\U_d} e^{x\Tr(U+U^*)}\det(U^*)^q dU,
\end{equation}
which we deduce from the following remarkable result of Wei and Wettig.

\begin{theorem}{\cite{w-w}}
\label{CFT}
For any matrices $X,Y \in \M_{(d+q) \times d}$ verifying $\det(Y^*X) \neq 0$ the
following holds:
\begin{equation}
	\int_{\U_{d+q}}e^{\Tr(Y^*UX+X^*U^*Y)}dU = H_{d \times q}
	\int_{\U_d} e^{\Tr(UY^*Y+X^*XU^*)}\det(UY^*X)^{-q} dU.
\end{equation}
\end{theorem}

To get (\ref{identity}) from Theorem \ref{CFT}, choose
\begin{equation}
	X=Y=\begin{bmatrix}
		t & 0 & \dots & 0\\
		0 & t & \dots & 0 \\
		0 & 0 & \dots & t \\
		0 & 0 & \dots & 0 \\
		\vdots & \vdots & \dots & \vdots \\
		0 & 0 & \dots & 0
		\end{bmatrix}
\end{equation}
where $t$ is a non-zero real number.  Then $Y^*X=t^2I,$ and 
$Y^*UX=T^{(q)}(U)$ for any unitary matrix $U \in \U_{d+q}.$  Thus Theorem \ref{CFT}
reduces to 
\begin{equation}
	\label{non-zero}
	\frac{t^{2dq}}{H_{d \times q}}\int_{\B_d} e^{t^2\Tr(P+P^*)}\gamma_d^{(q)}(dP) 
	=G_d(t^2;q) =  \int_{\U_d} e^{t^2\Tr(U+U^*)}\det(U^*)^q dU,
\end{equation}
valid for all real $t \neq 0.$
This identity also holds true at $t=0,$ since then the left hand side is obviously equal to 
$0$, and the right hand side is also equal $0$ for the following reason:
\begin{equation}
	\label{fact}
	\int_{\U_d} F(U)G(U^*) dU = 0
\end{equation}
for any two homogeneous polynomials $F,G$ in the entries of $U$ and $U^*$ with
$\deg F \neq \deg G$ (this is a standard property of Haar measure, see e.g. \cite{c-s} for
a proof).  Thus (\ref{non-zero}) is true for all real numbers $t$ and we obtain 
(\ref{identity}).

To complete the proof of Theorem \ref{main}, we compute the power series expansion of 
the unitary matrix integral (\ref{identity}) using techniques from symmetric function theory
(see \cite{stanley:ec2}, Chapter 7).  Recall that a monotone walk
\begin{equation}
	\emptyset \nearrow \nearrow \dots \nearrow \lambda
\end{equation}
on the Young graph $\Y$ is called a {\it standard Young tableau} of shape $\lambda.$
Following \cite{stanley:ec2}, we will use the notation
\begin{equation}
	f^{\lambda} := Z_{\Y}(R^{|\lambda|};\emptyset, \lambda)
\end{equation} 
for the number of standard Young tableaux of shape $\lambda.$

Expanding the integral (\ref{identity}) as a power series in $x,$ we have
\begin{equation}
	G_d(x;q)=\sum_{m,n \geq 0} \int_{\U_d} (\Tr U)^m (\Tr U^*)^n \det(U^*)^q dU
	\frac{x^{m+n}}{m!n!}.
\end{equation}
Appealing again to (\ref{fact}), this becomes
\begin{equation}
	G_d(x;q)=\sum_{n \geq 0}I_d(n;q)\frac{x^{2n+dq}}{n!(n+dq)!},
\end{equation}
where
\begin{equation}
	I_d(n;q) := \int_{\U_d} (\Tr U)^{n+dq}(\Tr U^*)^n \det(U^*)^qdU.
\end{equation}

Now we evaluate the integral $I_d(n;q).$
Let $\Lambda_d$ be the $\C$-algebra of symmetric polynomials
in $d$ indeterminates, and define the {\it Hall scalar product} on $\Lambda_d$ by
\begin{equation}
	\langle f | g \rangle:=\int_{\U_d} f(U)g(U^*) dU
\end{equation}
where $f(U),g(U^*)$ denote symmetric polynomials $f,g \in \Lambda_d$ evaluated on the 
spectra of $U$ and $U^*$ respectively.  The Hall product is a symmetric, non-negative 
bilinear form on $\Lambda_d$ (see \cite{stanley:ec2}).  It is well-known that the Schur polynomials 
$\{s_{\lambda}\}_{\lambda \in \Y_d}$ constitute a linear basis of $\Lambda_d,$ and 
moreover since each map $U \mapsto s_{\lambda}(U)$ is the character of an irreducible
polynomial representation of $\U_d$ we have
\begin{equation}
	\langle s_{\lambda} | s_{\mu} \rangle = \delta_{\lambda,\mu}
\end{equation}
by Schur orthogonality (see e.g. \cite{bump}).

Now
\begin{equation}
	I_d(n;q)=\langle e_1^{n+dq} | e_1^n e_d^q \rangle,
\end{equation}
where $e_1,e_d \in \Lambda_d$ are the first and last elementary symmetric polynomials.
In terms of Schur functions, on has the linear expansion
\begin{equation}
	e_1^N=\sum_{\substack{\lambda \in \Y_d\\ |\lambda|=n}}f^{\lambda}s_{\lambda}
\end{equation}
for powers of the first elementary symmetric polynomial (see \cite{stanley:ec2}).  It is also 
clear from the combinatorial definition of Schur polynomials in terms of semistandard Young
tableaux that $e_d^q=s_{d \times q}.$  Hence we have
\begin{equation}
	I_d(n;q) = \sum_{\substack{\mu, \lambda \in \Y_d\\ |\mu|=n+dq,|\lambda|=n}}
	f^{\mu}f^{\lambda} \langle s_{\mu}| s_{\lambda}s_{d \times q} \rangle.
\end{equation}

The following property of the Hall scalar product is well-known:
\begin{equation}
	\langle s_{\mu} |s_{\lambda} s_{\nu} \rangle = \langle s_{\mu/\nu} | s_{\lambda} \rangle,
\end{equation}
where $s_{\mu/\nu}$ is a skew Schur polynomial, which is by definition the zero polynomial
unless $\mu \supseteq \nu.$  In other words, the adjoint of the ``multiplication by 
$s_{\nu}$'' operator is the ``deletion of $\nu$'' operator.  Hence we have
\begin{equation}
	I_d(n;q) = \sum_{\substack{\mu,\lambda \in \Y_d \\ |\mu|=n+dq, |\lambda|=n}}
	f^{\mu}f^{\lambda} \langle s_{\mu/d \times q} | s_{\lambda} \rangle
	=\sum_{\substack{\lambda \in \Y_d \\ |\lambda|=n}}f^{\lambda+d \times q}f^{\lambda},
\end{equation}
where $\lambda+d \times q$ is the concatenation of the Young diagram $\lambda$ with 
the rectangular diagram $d \times q.$  Thus
\begin{align}
	I_d(n;q) &= \sum_{\substack{\lambda \in \Y_d \\ |\lambda|=n}} 
	Z_{\Y_d}(R^{n+dq};\emptyset,\lambda+d\times q)Z_{\Y_d}(R^n;\emptyset,\lambda) \\
	&=\sum_{\substack{\lambda \in \Y_d \\ |\lambda|=n}}
	Z_{\Y_d}(R^{n+dq};\emptyset,\lambda+d \times q)Z_{\Y_d}(L^n;\lambda+d\times q,d \times q) \\
	&=Z_{\Y_d}(L^nR^{n+dq};\emptyset,d \times q) \\
	&={2n +dq \choose n}^{-1}Z_d(2n+dq;q),
\end{align}
and we conclude that
\begin{align}
	G_d(x;q) &= \sum_{n \geq 0} I_d(n;q) \frac{x^{2n+dq}}{(n+dq)!n!} \\
	&= \sum_{n \geq 0} {2n +dq \choose n}^{-1}Z_d(2n+dq;q) \frac{x^{2n+dq}}{(n+dq)!n!} \\
	&=  \sum_{n \geq 0} Z_d(2n+dq;q) \frac{x^{2n+dq}}{(2n+dq)!} \\
	&= \sum_{N \geq 0} Z_d(N;q) \frac{x^N}{N!}.
\end{align}

\section{Asymptotics}
\label{asymp}
Let us now extract an explicit formula for $Z_d(N;q)$ in terms of random contractions from
Theorem \ref{main}.  We have
\begin{align}
	G_d(x;q) &= \frac{x^{dq}}{H_{d \times q}} \int_{\B_d} e^{x\Tr(P+P^*)} \gamma_d^{(q)}(dP) \\
	&= \sum_{N \geq 0} \frac{x^{dq}}{H_{d\times q}} 
	\int_{\B_d} (\Tr(P+P^*))^N \gamma_d^{(q)}(dP) \frac{x^N}{N!} \\
	&= \sum_{n \geq 0} \frac{{2n \choose n}}{H_{d \times q}}
	\int_{\B_d} |\Tr P|^{2n} \gamma_d^{(q)}(dP) \frac{x^{2n+dq}}{(2n)!}, \text{ by } (\ref{fact})\\
	&= \sum_{n \geq 0} \frac{(2n+dq)!}{H_{d \times q}n!n!} \int_{\B_d} |\Tr P|^{2n} 
	\gamma_d^{(q)}(dP) \frac{x^{2n+dq}}{(2n+dq)!}.
\end{align}
Thus 
\begin{equation}
	\label{explicit}
	Z_d(N;q) = \bigg{[} \frac{x^N}{N!} \bigg{]} G_d(x;q) = 
	\begin{cases}
		\frac{(2n+dq)!}{H_{d \times q}n!n!} \int_{\B_d} |\Tr P|^{2n} \gamma_d^{(q)}(dP), 
		\text{ if $N=2n+dq$}\\
		0, \text{ otherwise}
	\end{cases}.
\end{equation}
From this expression, we can easily obtain an asymptotic form for $Z_d(2n+dq;q)$ in the 
large $q$ limit with $d,n$ fixed.  This corresponds to the ``walk to infinity'' with bounded backtracking
Fisher's model.

Recall the following classical result of E. Borel \cite{borel}.  Let $p=(p_1,\dots,p_d,p_{d+1},\dots,
p_{d+q})$ be a uniformly random point from the real sphere $\mathbb{S}^{d+q-1}
\subset \R^{d+q}.$  Then, as $q \rightarrow \infty,$ 
$\sqrt{q}p_1,\dots,\sqrt{q}p_d$ converge weakly to an independent family of standard 
Gaussian random variables.  There is an ample generalization of this result to the setting of
truncated random unitary matrices, first established by Petz and R\'effy in \cite{p-r1}.

\begin{theorem}
\label{Borel}
Let $P_d^{(q)}$ be a random matrix from $\CUE^{(q)}.$  As $q \rightarrow \infty,$ 
the rescaled random matrix $\sqrt{q}P_d^{(q)}$ converges to the random matrix
$\Gamma_d$ whose entries are i.i.d. standard complex Gaussians (in the sense of 
pointwise convergence of spectral correlation functions).
\end{theorem}

We remark that this result of Petz and R\'effy follows readily from the 
work of Sommers and Zyczkowski \cite{s-z}.  In particular, one knows from
\cite{s-z} the the spectrum of a random matrix $P_d^{(q)}$ is a determinantal point
process in the unit disc $\mathbb{D}=\{z \in \C: |z| \leq 1\}$ governed by the 
Sommers-Zyczkowsi kernel
\begin{equation}
	\SZ_d^{(q)}(z,w) = \frac{q}{\pi} (1-|z|^2)^{\frac{q-1}{2}}(1-|w|^2)^{\frac{q-1}{2}}
	\sum_{j=0}^{d-1} {q+j \choose j}z^j \overline{w}^j,
\end{equation}
for any $q \geq 1.$  One then observes that the re-scaled kernel 
$q^{-1}SZ_d^{(q)}(q^{-\frac{1}{2}}z,q^{-\frac{1}{2}}w)$ converges to the Ginibre
kernel
\begin{equation}
	\operatorname{Gin}_d(z,w)
	= \frac{1}{\pi}e^{-(|z|^2+|w|^2)} \sum_{j=0}^{d-1}\frac{1}{j!}z^j \overline{w}^j
\end{equation}
as $q \rightarrow \infty.$  An elementary geometric proof of Theorem \ref{Borel} can
be found in \cite{m-t}, while a stronger assertion has recently been proved by
Krishnapur \cite{krishnapur}.

Theorem \ref{Borel} immediately implies that
\begin{equation}
	\int_{\B_d} |\Tr P|^{2n} \gamma_d^{(q)}(dP) \sim \frac{\E[|z_1+\dots+z_d|^{2n}]}{q^n}
\end{equation}
as $q \rightarrow \infty,$ where $z_1,\dots,z_n$ are independent standard complex Gaussians
and $\E$ denotes expected value.  Since the moment sequence of a standard complex Gaussian $z$
is well-known to be 
\begin{equation}
	\E[z^m \overline{z}^n] = \delta_{m,n}n!,
\end{equation}
and $z_1+\dots+z_d$ is a complex Gaussian of variance $d,$ we have
\begin{equation}
	\int_{\B_d} |\Tr P|^{2n} \gamma_d^{(q)}(dP) \sim \frac{d^nn!}{q^n}
\end{equation}
as $q \rightarrow \infty.$  Applying this asymptotic form in (\ref{explicit}), it follows from 
Stirling's approximation that
\begin{equation}
	\label{asymptotic}
	Z_d(N;q) \sim (2\pi)^{\frac{1-d}{2}} \bigg{(} \prod_{i=0}^{d-1} i! \bigg{)} 
	d^{3n+dq+\frac{1}{2}}q^{n+\frac{1-d^2}{2}}
\end{equation}
as $q \rightarrow \infty,$ where $N=2n+dq$ and $n \geq 0, d \geq 1$ are fixed
but arbitrary.  This asymptotic form reveals an interesting
phase transition in the large $q$ asymptotic behaviour of $Z_d(N;q)$ at $n \approx \frac{d^2-1}{2}$
from 
\begin{equation}
	\frac{\text{exponential in $q$}}{\text{polynomial in $q$}}
\end{equation}
to
\begin{equation}
	(\text{exponential in $q$})(\text{polynomial in $q$}).
\end{equation}

\end{document}